\RequirePackage{fix-cm}
\documentclass[smallextended]{svjour3}       
\smartqed  
\usepackage{graphicx}
\usepackage{amsmath,amssymb,amsfonts}
\usepackage{dsfont} 
\usepackage[colorlinks=true,linkcolor=blue,citecolor=blue,urlcolor=blue]{hyperref}

\newcommand{\dd}{\mathrm{d}}
\newcommand{\1}{\mathds{1}}

\numberwithin{equation}{section}

\begin{document}

\title{Dependence Properties of B-Spline Copulas \\
 \thanks{This work was supported by JSPS KAKENHI Grants, Numbers 16K00060 and 16H02792.}
}
\subtitle{Dependence properties of B-spline copulas}


\author{Xiaoling Dou \and
        Satoshi Kuriki \and \\
        Gwo Dong Lin \and
        Donald Richards 
}


\institute{X. Dou \at
              Waseda University, 3-4-1 Ohkubo, Shinjuku, Tokyo 169-8555, Japan \\
              Tel.: +81-3-5286-2358\\
              \email{xiaoling@aoni.waseda.jp}           
           \and
           S. Kuriki \at
              The Institute of Statistical Mathematics, 10-3 Midoricho, Tachikawa, Tokyo 190-8562, Japan
          \and
          G. D. Lin \at
          Institute of Statistical Science, Academia Sinica, Taipei 11529, Taiwan, R.O.C.
          \and
          D. Richards \at
          Department of Statistics, Pennsylvania State University, University Park, PA 16802, U.S.A.
}

\date{Received: date / Accepted: date}

\maketitle

\begin{abstract}
We construct by using B-spline functions a class of copulas that
includes the Bernstein copulas arising in Baker's distributions.  
The range of correlation of the B-spline copulas is examined, and 
the Fr\'{e}chet--Hoeffding upper bound is proved to be attained when 
the number of B-spline functions goes to infinity.  As the B-spline 
functions are well-known to be an order-complete weak Tchebycheff
system from which the property of total positivity of any order
(TP${}_\infty$) follows for the maximum correlation case, the 
results given here extend classical results for the Bernstein
copulas.  In addition, we derive in terms of the Stirling numbers 
of the second kind an explicit formula for the 
moments of the related B-spline functions on $[0,\infty)$.  
\keywords{
Bernstein copula 
\and Fr\'{e}chet--Hoeffding upper bound 
\and Order-complete weak Tchebycheff system 
\and Schur function 
\and Stirling number of the second kind
\and Total positivity of order $r$}
\end{abstract}

\section{Introduction: A review of the Bernstein copulas}

A novel method that applied the theory of order statistics to construct
multivariate distributions with given marginal distributions was
developed by Baker \cite{Baker08}.  We refer to Lin, et al. \cite{Lin-etal14} 
for a recent survey of this topic.  

Baker's idea, applied to univariate 
cumulative distribution functions $F$ and $G$, can be described as
follows: Let $\{X_1,\ldots,X_n\}$ and $\{Y_1,\ldots,Y_n\}$ be independent 
random samples from the distributions $F$ and $G$, respectively. Let 
$X_{k,n}$ be the $k$th smallest order
statistic of $\{X_1,\ldots,X_n\}$, and denote by $F_{k,n}$ the distribution 
of $X_{k,n}$; we write this as $X_{k,n}\sim F_{k,n}$.  Similarly, we denote 
by $Y_{k,n}$ the $k$th smallest order statistic of $\{Y_1,\ldots,Y_n\}$ and 
we let $G_{k,n}$ be its corresponding distribution, written 
$Y_{k,n}\sim G_{k,n}$.  (Note that $F$ and $G$ can be discrete 
distributions.) 

Let $R = (r_{k\ell})_{1\le k,\ell\le n}$ be a parameter matrix whose 
matrix entries $r_{k\ell}$ satisfy the conditions 
\begin{equation}
\label{R0}
\sum_{k=1}^n r_{k\ell} = \sum_{\ell=1}^n r_{k\ell} = \frac{1}{n}, \quad
 r_{k\ell} \ge 0,\ k,\ell=1,2,\ldots,n.
\end{equation}
Now choose the pair $(X_{k,n},Y_{\ell,n})$ with probability $r_{k\ell}$,
$k,\ell=1,2,\ldots,n$. Then $(X_{k,n},Y_{\ell,n})$ 
follows Baker's bivariate  distribution: For $x,y\in\mathbb{R}$, the 
joint cumulative distribution function 
$H(x,y) := \Pr(X_{k,n} \le x,\,Y_{k,n} \le y)$ satisfies 
\begin{equation}
\begin{aligned}
\label{H}
 H(x,y)
 &= \sum_{k=1}^{n}\sum_{\ell=1}^{n} r_{k\ell} F_{k,n}(x) G_{\ell,n}(y) \\
 &= (F_{1,n}(x),\ldots,F_{n,n}(x))\,R\,
 (G_{1,n}(y),\ldots,G_{n,n}(y))^{\!\mathsf{T}}, 
\end{aligned}
\end{equation}
where ``${\mathsf{T}}$'' denotes transpose.  It is immediately evident that 
$H$ has marginal distributions $F$ and $G$.

Let $B_{k,n}$ be the distribution function of the $k$th
smallest order statistic of a random sample of size $n$
from the uniform distribution $\mathcal{U}$ on $[0,1]$.  
It is well-known that 
\begin{equation}
\label{b-density}
 B_{k,n}(u) = \int_0^u b_{k,n}(t)\,\dd t
 \end{equation}
where
$$
b_{k,n}(t)= n \binom{n-1}{k-1} t^{k-1} (1-t)^{n-k},
$$
$t \in [0,1]$.  Furthermore, $F_{k,n}$ equals the composition 
$B_{k,n}\circ F$ (see, e.g., Hwang and Lin \cite{Hwang-Lin84})
and Baker's bivariate  distribution (\ref{H}) 
can be rewritten as
\[
 H(x,y) = C(F(x),G(y);R),
\]
$x, y \in \mathbb{R}$, where, for $u, v \in [0,1]$, 
\begin{equation}
\label{C0}
 C(u,v;R) = \sum_{k=1}^{n}\sum_{\ell=1}^{n} r_{k\ell} B_{k,n}(u) B_{\ell,n}(v)
\end{equation}
is a copula function with parameter matrix $R$ satisfying (\ref{R0}).  
Conversely, if the marginals $F$ and $G$ are equal to $\mathcal{U}$ then 
Baker's bivariate distribution (\ref{H}) reduces to the copula (\ref{C0}).

The copula $C(u,v;R)$ in (\ref{C0}) is called the \textit{Bernstein 
copula with parameter matrix} $R$ because $b_{k,n}/n$ is a Bernstein 
polynomial (see Dou, et al. \cite{Dou-etal14}).  By differentiating 
(\ref{C0}) with respect to $u$ and $v$, we obtain the Bernstein copula 
density: 
\begin{equation*}
 c(u,v;R) = \sum_{k=1}^{n}\sum_{\ell=1}^{n} r_{k\ell} b_{k,n}(u) b_{\ell,n}(v),
\end{equation*}
$u, v \in [0,1]$.

Within the parameter space (\ref{R0}) of $R$, the maximum correlation
is attained when $r_{k\ell}=(1/n)\delta_{k\ell}$, i.e., $R=(1/n)I_n$:
\[
 C^*(u,v) := C(u,v;(1/n)I_n) = \frac{1}{n} \sum_{k=1}^{n} B_{k,n}(u) B_{k,n}(v),
\]
with corresponding density
\begin{equation*}
 c^*(u,v) := c(u,v;(1/n)I_n) = \frac{1}{n} \sum_{k=1}^{n} b_{k,n}(u) b_{k,n}(v),
\end{equation*}
$u, v \in [0,1]$.

Dou, et al. \cite{Dou-etal13} proved that the maximum correlation copula 
$C^*(u,v)$ and its density $c^*(u,v)$ both are totally positive of order 
$2$ (TP${}_2$) in $(u,v)$ \cite{Karlin-Studden66,Karlin68,Pinkus10}.  One 
of the main purposes of the present paper is to show further that both 
$C^*$ and $c^*$ are TP${}_{\infty}$, i.e., TP${}_r$ for all $r\ge 2$.

In Section \ref{sec:b-spline}, we introduce first the general 
order-complete weak Tchebycheff (OCWT) systems 
 and then a class of copulas, based on B-spline functions,
 that includes the Bernstein copulas $C(u,v;R)$ in (\ref{C0}).
The maximum correlation copula $C^*(u,v)$ and its total positivity 
properties are investigated in Sections \ref{sec:maximum_correlation} 
and \ref{sec:total_positivity}, respectively.
Finally, in Section \ref{sec:moment} we calculate the moments of the
related B-spline functions on $[0,\infty)$ and 
make the connection with the Stirling numbers of the second kind.

\section{B-spline copulas}
\label{sec:b-spline}
\setcounter{equation}{0}

We consider first a general setting based on OCWT systems, and then we 
define a class of B-spline copulas that includes the Bernstein copulas 
as special cases.  In other words, we will show that a larger class of 
candidate copulas can play the same roles as $B_{k,n}$ and $R$ while 
still retaining the desired properties of the copula $C(u,v;R)$.

Let $q_k \ge 0$, $k=1,\ldots,n$, $\sum_{k=1}^n q_k=1$, and let
$\phi_1,\ldots,\phi_n$ be probability densities on $[0,1]$ such that
\begin{equation}
\label{unity}
 \sum_{k=1}^n q_k \, \phi_{k}(t) = 1,
\end{equation}
$t \in [0,1]$.  We assume further that $\{\phi_1,\ldots,\phi_n\}$ is an 
{\it order-complete weak Tchebycheff system} (OCWT-system), i.e., 
\begin{itemize}
\item[(i)] $\phi_1,\ldots,\phi_n$ are linearly independent, 
and 
\item[(ii)] $\phi_k(t)$ is \textit{totally positive of order $n$ (TP$_n$) 
in $(k,t)$}, i.e., for each $r = 1,\ldots,n$, 
\begin{equation}
\label{wT}
 \det\big(\phi_{k_i}(t_j)\big)_{r\times r} \ge 0
\end{equation}
for all $k_1>\cdots>k_r$ and $t_1>\cdots>t_r$.
\end{itemize}
See Karlin and Studden \cite[Chapter 1]{Karlin-Studden66} or Schumaker 
\cite[Chapter 2]{Schumaker07} for examples of OCWT systems.

Let $q_{1k} \ge 0$, $k=1,\ldots,n_1$, such that $\sum_{k=1}^{n_1}
q_{1k}=1$. Also, let $q_{2\ell}\ge 0$, $\ell=1,\ldots,n_2$, such that
$\sum_{\ell=1}^{n_2} q_{2\ell}=1$.
Letting
\[
 \Phi_{k}(u) = \int_0^u \phi_{k}(t)\,\dd t,
\]
$u \in [0,1]$, we define the \textit{B-spline copula}, a generalization of 
the Bernstein copula (\ref{C0}), by
\begin{equation}
\label{C}
 C(u,v;R) = \sum_{k=1}^{n_1} \sum_{\ell=1}^{n_2} r_{k\ell} \, \Phi_{k}(u) \, \Phi_{\ell}(v),
\end{equation}
$u, v \in [0,1]$, with parameter matrix
\begin{equation}
\begin{aligned}
& R=(r_{k\ell})_{1\le k\le n_1;\,1\le\ell\le n_2},\ \ r_{k\ell}\ge 0, \\
& \sum_{k=1}^{n_1} r_{k\ell}= q_{2\ell},\ \sum_{\ell=1}^{n_2} r_{k\ell}=q_{1k},\ \ 
k=1,2,\ldots,n_1,\ \ell=1,2,\ldots,n_2.
\label{R}
\end{aligned}
\end{equation}
The copula (\ref{C}) is a \textit{bona fide} copula since, for any $u\in[0,1]$, 
\begin{align*}
C(u,1;R)
 &= \sum_{k=1}^{n_1}\sum_{\ell=1}^{n_2} r_{k\ell} \Phi_{k}(u) = \sum_{k=1}^{n_1} q_{1k} \Phi_{k}(u) \\
 &= \int_0^{u} \sum_{k=1}^{n_1} q_{1k} \phi_{k}(t)\,\dd t = \int_0^{u} 1\,\dd t = u;
\end{align*}
and similarly, $C(1,v;R) = v$, $v\in[0,1]$.

Throughout the paper, we restrict our attention to the case in which $n_1=n_2 = n$ 
and $q_{1k}=q_{2k} = q_k$; further, we use the notation $Q = \mathrm{diag}(q_k)_{1\le k\le n}$
for the diagonal matrix with diagonal entries $q_1,\ldots,q_n$.

\begin{theorem}
\label{thm:maxcorr}
For the copula (\ref{C}) with the parameter
space (\ref{R}), the maximum correlation is attained when
$ r_{k\ell} = q_k \delta_{k\ell}$, equivalently, $R = Q$.
\end{theorem}

In the maximum correlation case, $C(u,v;R)$ becomes
\begin{equation}
\label{C*}
 C^*(u,v) := C(u,v;Q) = \sum_{k=1}^{n} q_k \Phi_{k}(u) \Phi_{k}(v),
\end{equation}
$u,v\in[0,1]$.  

To prove Theorem \ref{thm:maxcorr}, we need the following crucial
lemma.  This result is a generalization of the weak majorization inequality
on the closed simplicial cone
\[
 \mathcal{D}_+=\{(x_1,\ldots,x_n):x_1\ge \cdots\ge x_n\ge 0\}\subset\mathbb{R}^n_+
\]
for doubly stochastic matrices \cite[p. 639]{Marshall-Olkin-Arnold11}.

\begin{lemma}
\label{lem:majorization}
Let $a_1\ge\cdots\ge a_n\ge 0$ and
$b_1\ge\cdots\ge b_n\ge 0$ be given.  Let $q_1,\ldots,q_n\ge 0$
satisfy $\sum_{k=1}^nq_k=1$.  Then,
\begin{equation}
\label{lem}
\max_{\substack{\sum_k r_{k\ell}=q_\ell \\ \sum_\ell r_{k\ell}=q_k \\ r_{k\ell}\ge 0}} \, 
 \sum_{k=1}^n \sum_{\ell=1}^n r_{k\ell}a_k b_\ell\ =\ \sum_{k=1}^n q_ka_k b_k.
\end{equation}
\end{lemma}

\noindent
\textit{Proof.} 
Let  $\mathbf{a}=(a_1,\ldots,a_n)$ and $\mathbf{b}=(b_1,\ldots,b_n)$.
Define 
$$
p_{k\ell} = 
\begin{cases}
r_{k\ell}, & k\ne \ell , \\
r_{kk}+1-q_k, & k=l .
\end{cases}
$$
For given $Q = \mathrm{diag}(q_k)_{1\le k\le n}$, it is 
straightforward to see that $P = R+I_n-Q$ is a $n\times n$ doubly 
stochastic matrix.  Hence, by the famous characterization of majorization 
due to Hardy, Littlewood, and P\'olya \cite{HLW29}, the vector 
$\mathbf{c} = \mathbf{a}P$ is majorized by $\mathbf{a}$, denoted 
$\mathbf{c}\prec\mathbf{a}$.

We now rearrange the components of $\mathbf{c}=(c_1,\ldots,c_n)$
in decreasing order, listing them as $c_{[1]}\ge\cdots\ge c_{[n]}$,
and let $\mathbf{c}^*=(c_{[1]},\ldots,c_{[n]})$.
Then we have also $\mathbf{c}^*\prec\mathbf{a}$ and hence
$\mathbf{c}^*\mathbf{b}^{\!\mathsf{T}}\le\mathbf{a}\mathbf{b}^{\!\mathsf{T}}$
because $\mathbf{a}, \mathbf{b}, \mathbf{c}^*\in\mathcal{D}_+$
 (see \cite[p.\,133]{Marshall-Olkin-Arnold11}).
On the other hand, note that
$\mathbf{c}\mathbf{b}^{\!\mathsf{T}}\le\mathbf{c}^*\mathbf{b}^{\!\mathsf{T}}$
because $\mathbf{b}, \mathbf{c}^*\in\mathcal{D}_+$.
These results together imply that
$\mathbf{a}P\mathbf{b}^{\!\mathsf{T}}
=\mathbf{c}{\mathbf{b}^{\!\mathsf{T}}}\le\mathbf{c}^*\mathbf{b}^{\!\mathsf{T}}\le\mathbf{a}\mathbf{b}^{\!\mathsf{T}}
$.

Consequently, $\max_P\mathbf{a}P{\mathbf{b}^{\!\mathsf{T}}}= \mathbf{a}\mathbf{b}^{\!\textsf{T}}$, which we can write alternatively as 
\[
 \max_{\substack{\sum_k p_{k\ell}=1 \\ \sum_\ell p_{k\ell}=1 \\ p_{k\ell}\ge 0}} \
 \sum_{k=1}^n \sum_{\ell=1}^n p_{k\ell}a_k b_\ell\ =\ \sum_{k=1}^n a_k b_k.
\]
Equivalently, $\max_R\mathbf{a}R\mathbf{b}^{\!\mathsf{T}}= \mathbf{a}Q\mathbf{b}^{\!\mathsf{T}}$
by cancelling the common term $\mathbf{a}\mathbf{b}^{\!\mathsf{T}}$ on both sides above.  This 
is exactly (\ref{lem}), so the proof now is complete.
\hfill $\qed$

\medskip

\noindent
\textit{Proof of Theorem \ref{thm:maxcorr}.}  
Since $\{\phi_1,\ldots,\phi_n\}$ is an OCWT-system then, for all $i<j$
and $s<t$,
\[
 \phi_i(s) \phi_j(t) - \phi_j(s) \phi_i(t) = \det\begin{pmatrix} \phi_i(s) & \phi_i(t) \\ \phi_j(s) & \phi_j(t) \end{pmatrix} \ge 0.
\]
Integrating this inequality with respect to $(s,t)$ over $s\in (0,u)$ and $t\in (u,1)$, we obtain 
\[
 \Phi_i(u)(1-\Phi_j(u)) - \Phi_j(u)(1-\Phi_i(u)) =
 \Phi_i(u) - \Phi_j(u) \ge 0,
\]
$u\in[0,1]$.  Therefore, we obtain the stochastic order, 
\[
 \Phi_{1}(u) \ge \Phi_{2}(u) \ge \cdots \ge \Phi_k(u)
\]
for all $u\in[0,1]$.  Combining this result with Lemma \ref{lem:majorization}, we obtain
the inequality
\[
 C^*(u,v) := C(u,v;Q) \ge C(u,v;R)
\]
for all $u,v\in[0,1]$ and $R$ satisfying (\ref{R}).  The theorem now follows from Hoeffding's covariance formula, 
\[
 \mathrm{Cov}(X,Y) = \int_{-\infty}^{\infty}\int_{-\infty}^{\infty}\big[\Pr(X\le x,Y\le y)-\Pr(X\le x)\Pr(Y\le y)\big]\,\dd x \dd y
\]
(see, e.g., \cite{Lin-etal14}).  The proof is complete.
\hfill $\qed$
\bigskip

Functions $\phi_k$ satisfying (\ref{unity}) and (\ref{wT}) can be
constructed by B-spline functions as we now show.  Let $N^d_i$ be a
B-spline function on $[0,1]$ of degree $d\ (\ge 0)$ defined as a
non-zero B-spline basis with $m+2d+1$ knots:
\begin{equation}
\label{knots}
 \underbrace{t_{-d} = \cdots = t_{-1}}_d = t_0 = 0 < t_1 < \cdots 
 < t_{m-1} < 1 = t_m = \underbrace{t_{m+1} = \cdots = t_{m+d}}_d.
\end{equation}
Then, $N^d_{i}(t)$ is generated by the recursion formula, 
\begin{align*}
 N^d_i(t) &= \frac{t-t_i}{t_{i+d}-t_i}N^{d-1}_i(t) + \frac{t_{i+d+1}-t}{t_{i+d+1}-t_{i+1}}N^{d-1}_{i+1}(t) \\
 &= \frac{t-t_i}{t_{i+d}-t_i}N^{d-1}_i(t) + \biggl(1-\frac{t- t_{i+1}}{t_{i+d+1}-t_{i+1}}\biggr)N^{d-1}_{i+1}(t),
\end{align*}
$t\in[0,1]$, for $i=-d,\ldots,-1,0,1,\ldots,m-1$, with initial conditions
\[
 N^0_i(t) = \begin{cases} \displaystyle
 1, & i<m \ \mbox{and}\ t\in [t_i,t_{i+1}), \\
    & \mbox{or}\ i=m-1 \ \mbox{and}\ t=t_m=1, \\
 0, & \mbox{otherwise} \end{cases}
\]
(see \cite{deBoor72,deBoor01,Nurnberger89}).
The number of non-zero bases is
\[
 n=m+d.
\]

The B-spline is known to satisfy 
\begin{itemize}
\item[(i)] $N^d_i(t)\ge 0,\ t\in[0,1]$,
\item[(ii)] The support is given by
\[
 \mathop{\mathrm{supp}}N^d_i = \overline{\{ t \mid  N^d_i(t)>0 \}} = [t_i,t_{i+d+1}],
\]
$i=-d,\ldots,-1,0,1,\ldots,m-1$, and 
\item[(iii)] The ``partition of unity'' property: 
\[
 \sum_{i=-d}^{m-1} N^d_i(t) = 1\ \ \mbox{for all}\ t\in [0,1].
\]
\end{itemize}
For given $d$ and $m$, let
\begin{equation}
\label{B-spline}
q_{k} =q_{k,d} = \int_0^1 N^d_{k-d-1}(t)\,\dd t, \quad \phi_{k}(t)=\phi_{k,d}(t) = \frac{1}{q_{k}}
N^d_{k-d-1}(t),
\end{equation}
where $t\in[0,1]$ and $k=1,2,\ldots,n\ (=m+d)$.  Then, (\ref{unity})
holds, and we have the following result (see \cite{deBoor76}, or
\cite[Theorems 4.18 and 4.65]{Schumaker07}).

\begin{theorem}
\label{thm:WT} Under the hypotheses (\ref{knots}) and (\ref{B-spline}), 
the set 
$\{N^d_i\}_{i=-d}^{m-1}$ of B-spline
functions, and hence also the B-spline system $\{\phi_1,\ldots,\phi_n\}$, forms
an OCWT-system satisfying (\ref{wT}).
\end{theorem}

To illustrate the use of B-spline systems, we now provide some examples.

\begin{theorem}
\label{thm:m=1}
 Let $m=1$ and the degree $d=n-1 (= n-m)$.  Then the B-splines (\ref{B-spline}) reduce to the
Bernstein system (\ref{b-density}).  Specifically, for $k=1,\ldots,n$ and $t\in[0,1]$, 
\begin{equation}
\label{bernstein}
 q_k =q_{k,d}= \frac{1}{n}, \quad \phi_{k}(t) =\phi_{k,d}(t)= b_{k,n}(t).
\end{equation}
\end{theorem}

\noindent
\textit{Proof.}
We prove the result by induction on $d$.
Note that for $d=0$ (i.e., $n=1$), $N_0^0(t)=1$, $t\in[0,1]$, and hence
\[
 q_1=q_{1,0}=\int_0^1 N_0^0(t)\,\dd t=1, \quad
 \phi_1(t)=\phi_{1,0}(t)=N_0^0(t)=b_{1,1}(t),
\]
$t\in[0,1]$.  For $d=1$ (i.e., $n=2$), we have the required $N_i^1$, $q_{k,1}$, and
$\phi_{k,1}$ as follows:
\[
 N_{-1}^1(t)=(1-t)N_0^0(t)=1-t, \quad N_0^1(t)=tN_0^0(t)=t,\ \ t\in[0,1],
\]
\[
 q_1=q_{1,1}=\int_0^1 N_{-1}^1(t)\,\dd t=1/2, \quad q_2=q_{2,1}=\int_0^1 N_{0}^1(t)\,\dd t=1/2,
\]
\[
 \phi_{1}(t) =\phi_{1,1}(t)=2N_{-1}^1(t)= b_{1,2}(t),
\]
\[
 \phi_{2}(t) =\phi_{2,1}(t)=2N_{0}^1(t)= b_{2,2}(t),\ \ t\in[0,1].
\]
Assume that the theorem holds true for $d=n-2$, then we  want to
prove (\ref{bernstein}) for $d=n-1$.  In this case, the B-spline
functions are of the form
\begin{align*}
N_{k-n}^{n-1}(t)
& = \begin{cases}
\displaystyle
  (1-t)N_{2-n}^{n-2}(t)={(n-1)}^{-1}(1-t)b_{1,n-1}(t),& k=1, \\[2mm]
 tN_{k-n}^{n-2}(t)+(1-t)N_{k-n+1}^{n-2}(t) \\
 \quad ={(n-1)}^{-1}[tb_{k-1,n-1}(t)+(1-t)b_{k,n-1}(t)],& 1<k<n, \\[2mm]
  tN_0^{n-2}(t)={(n-1)}^{-1}tb_{n-1,n-1}(t), & k=n.
\end{cases}
\end{align*}
It can be shown that for $k=1,2,\ldots,n$,
$$
q_k=q_{k,n-1}=\int_0^1 N_{k-n}^{n-1}(t)\,\dd t=1/n
$$
and 
$$
N_{k-n}^{n-1}(t)=q_kb_{k,n}(t),
$$
$t\in[0,1]$.  This completes the proof.
\hfill $ \qed$

\bigskip

From now on, for simplicity, we consider only the B-spline with
equally-spaced knots, i.e., the B-spline functions on $[0,1]$ of
order $d$ having knots given in (\ref{knots}) with $t_i = i/m$, 
$i=1,2,\ldots,m-1$.

\begin{example}
Suppose $d=0$, i.e., $n=m$; then the B-spline system becomes a ``histogram''.
Namely, for $k=1,2,\ldots,n$, 
\begin{equation}
\label{d=0}
 q_k =q_{k,0}= \frac{1}{n}, \quad \phi_{k}(t) =\phi_{k,0}(t)=n\1_{[\frac{k-1}{n},\frac{k}{n})}(t),
\end{equation}
$t \in [0,1]$, where $\1_A$ denotes the indicator function of the set $A$.
\end{example}

\begin{example}
For $d=1$, i.e., $n=m+1$, we have 
$$
q_1=q_{1,1}= \frac{1}{2m} \, , \ \ \ \ 
q_2=q_3=\cdots=q_{m-1}=q_m=\frac{1}{m}\, , \ \ \ \ 
q_n=q_{n,1}=\frac{1}{2m}\, , 
$$
and the B-spline system is 
\begin{align*}
 \phi_1(t) &= \phi_{1,1}(t)=q_1^{-1}N_{-1}^1(t)=2m(1-mt)\1_{[0,\frac{1}{m})}(t), \\
 \phi_n(t) &= \phi_{n,1}(t)=q_n^{-1}N_{m-1}^1(t)=2m(mt-m+1)\1_{[1-\frac{1}{m},1)}(t),
\end{align*}
and, for $k=2,3,\ldots,n-1$,
\begin{align*}
 \phi_k(t) = \phi_{k,1}(t) &= q_k^{-1}N_{k-2}^1(t) \\
 &=q_k^{-1}[(mt-k+2)N_{k-2}^0(t)+(k-mt)N_{k-1}^0(t)] \\
 &= m\bigl[(mt-k+2)\1_{[\frac{k-2}{m},\frac{k-1}{m})}(t)+(k-mt)\1_{[\frac{k-1}{m},\frac{k}{m})}(t)\bigr],
\end{align*}
$t\in[0,1]$.  

We remark that Shen, et al. \cite{Shen-etal08} earlier proposed the ``linear
B-spline copula'', which corresponds to the case $d=1$.
\end{example}

\begin{example}
\label{ex:(n,d)=(5,3)}
For $d=3$ and $m=2$, i.e., $n=5$, we have 
\[
 q_1=1/8,\ \ q_2=q_3=q_4=1/4,\ \ q_5=1/8,
\]
and the B-spline system is 
\begin{align*}
 \phi_1(t) &= q_1^{-1}N_{-3}^3(t)=8(1-2t)^3\1_{[0,\frac{1}{2})}(t), \\
 \phi_2(t) &= q_2^{-1}N_{-2}^3(t)=8t(7t^2-9t+3)\1_{[0,\frac{1}{2})}(t) +8(1-t)^3\1_{[\frac{1}{2},1)}(t), \\
 \phi_3(t) &= q_3^{-1}N_{-1}^3(t)=8t^2(3-4t)\1_{[0,\frac{1}{2})}(t) +8(1-t)^2(4t-1)\1_{[\frac{1}{2},1)}(t), \\
 \phi_4(t) &= \phi_2(1-t), \\ 
 \phi_5(t) &= \phi_1(1-t),
\end{align*}
$t\in[0,1]$.  The means of the densities $\phi_1,\ldots,\phi_5$ are $1/10$,
$3/10$, $1/2$, $7/10$, and $9/10$, respectively.  We use these values in the 
computation of Table \ref{tab:rankcorr} of maximum correlations for $(n,d)=(5,3)$.
\end{example}

\section{The maximum correlation copula: Range of correlation}
\label{sec:maximum_correlation}
\setcounter{equation}{0}

For copula functions, the range of the correlation is of particular importance.  
In particular, great attention is paid to the maximum achievable correlation 
(see, e.g., Lin and Huang \cite{Lin-Huang10}).  By Theorem \ref{thm:maxcorr}, 
the maximum is attained when the copula density is
\begin{equation}
\label{c*}
 c^*(u,v) = \sum_{k=1}^{n} q_k \phi_{k}(u) \phi_{k}(v),\ \ u,v\in[0,1].
\end{equation}

Suppose that $(U,V)$ is from the copula density (\ref{c*}).  Then,
\begin{equation}
\label{pm}
 E[U V] = \sum_{k=1}^{n} q_k \biggl(\int_0^1 u \phi_{k}(u)\, \dd u\biggr)^2.
\end{equation}
Noting that $E[U]=E[V]=1/2$ and $\mathrm{Var}(U)=\mathrm{Var}(V)=1/12$, it follows that 
\begin{equation}
\label{corr}
 \mathrm{corr}(U,V) = 12 \Big(E[U V] - \frac{1}{4}\Big).
\end{equation}
In the Bernstein case ($m=1$), it follows from Theorem \ref{thm:m=1} that 
$E[UV]=(2n+1)/6(n+1)$ and hence
$$
\mathrm{corr}(U,V)=1-\frac{2}{n+1}.
$$
In the case of the B-spline of order zero, given in (\ref{d=0}),
\begin{equation*}
 \mathrm{corr}(U,V)
 = 12 \Biggl( \frac{1}{n}\sum_{k=1}^{n} \biggl(n \int_{(k-1)/n}^{k/n} t\,\dd t \biggr)^2 - \frac{1}{4}\Biggr) = 1 - \frac{1}{n^2}.
\end{equation*}
In particular, when $d=0$ and $n=m=1$, $\phi_1(t)=1$ on $[0,1]$, and
hence $C^*(u,v)=0,\ u,v\in[0,1]$.
This is the independent case, so $\mathrm{corr}(U,V)=0$.

In order to calculate the maximum correlation for general $d$, we
present first a lemma in which it is understood that
the vectors $(q_k)$ and $(r_k)$ reduce to the central parts when $d=0$.

\begin{lemma}
\label{lem:m}
Suppose that $m \ge d \ge 0$, i.e., $n=m+d \ge 2d \ge 0$.  Let
$N^d_i$, $i=-d,-d+1,\ldots,m-1$, be the B-spline functions on
$[0,1]$ of order $d$ having knots (\ref{knots}) with 
$t_i = i/m$, $i=0,1,\ldots,m$.  
In addition, denote the integral and the first moment of
$N^d_{k-d-1}$ by
\begin{equation*}
 q_k = \int_0^1 N^d_{k-d-1}(t)\,\dd t \ \ \ \ \hbox{and}\ \ \ \ 
 r_k = \int_0^1 t N^d_{k-d-1}(t)\,\dd t,
\end{equation*}
$k=1,\ldots,n$.  Then, 
\begin{align*}
 (q_k)_{1\le k\le n} =& \frac{1}{m}\biggl(\,\underbrace{\frac{1}{d+1},\frac{2}{d+1},\ldots,
 \frac{d}{d+1}}_d,\underbrace{1,\ldots,1}_{m-d},\underbrace{\frac{d}{d+1},\frac{d-1}{d+1},\ldots,\frac{1}{d+1}}_d\,\biggr), \\
 (r_k)_{1\le k\le n} =& \frac{1}{m^2}\biggl(\,\underbrace{\frac{1^2(1+1)}{2(d+1)(d+2)},\frac{2^2(2+1)}{2(d+1)(d+2)},\ldots,\frac{d^2(d+1)}{2(d+1)(d+2)}}_d, \\
& \qquad \underbrace{\frac{d+1}{2},\frac{d+3}{2},\ldots,\frac{2m-1-d}{2}}_{m-d},
\underbrace{m^2(q_d-r_d),\ldots,m^2(q_1-r_1)}_d\biggr).
\end{align*}
\end{lemma}

\noindent
\textit{Proof.}  
For $1\le k\le m$, we have $q_k=\gamma^d_{k-d-1}(0)/m$ and 
$r_k= \gamma^d_{k-d-1}(1)/m^2$, where $\gamma^d_i(0)$ and 
$\gamma^d_i(1)$ are given below in (\ref{h0}) and (\ref{h1}), 
respectively.  Also, for $k=m+1,\ldots,m+d = n$, we 
have the relations
\begin{equation}
\label{symmetry}
 q_k = q_{n+1-k}, \quad r_k = q_{n+1-k}-r_{n+1-k}\,,
\end{equation}
 because
$N_{k-d-1}^d(t)=N_{m-k}^d(1-t)=N_{(n+1-k)-d-1}^d(1-t)$, 
$t\in[0,1]$.  Solving the equations (\ref{symmetry}) in a 
successive manner, we obtain the stated results.  
\hfill $\qed$

\begin{theorem}
Under the assumptions of Lemma \ref{lem:m}, suppose that $(U,V)$ 
have the copula density $c^*$ in (\ref{c*}) with $\phi_k$ defined 
through the B-spline functions (\ref{B-spline}) having knots 
given in Lemma \ref{lem:m}. Then the correlation of $(U,V)$ is
\begin{align}
 \mathrm{corr}(U,V)
= 1 - \frac{d+1}{(n-d)^2} + \frac{d (d+3) (2 d+3)}{5 (d+2) (n-d)^3}.
\label{maxcorr}
\end{align}
\end{theorem}

\noindent
\textit{Proof.}  Using (\ref{pm}), (\ref{symmetry}), and the 
notations in Lemma \ref{lem:m}, write first
\begin{align}
\label{pm2}
E[UV]
= \sum_{k=1}^n \frac{r_k^2}{q_k}
&= \sum_{k=1}^m \frac{r_k^2}{q_k} + \sum_{k=1}^d \frac{(q_k-r_k)^2}{q_k}
 \nonumber\\
&= \sum_{k=1}^m \frac{r_k^2}{q_k} + \sum_{k=1}^d \Big(q_k - 2r_k
 + \frac{r_k^2}{q_k}\Big) \nonumber\\
&= 2\sum_{k=1}^d \frac{r_k^2}{q_k}
 + \sum_{k=1}^dq_k-2\sum_{k=1}^dr_k + \sum_{k=d+1}^m\frac{r_k^2}{q_k}.
\end{align}
The final result is obtained by substituting (\ref{pm2}) in (\ref{corr}) 
and carry out the calculations to obtain (\ref{maxcorr}) with the help of 
Lemma \ref{lem:m}.
\hfill $\qed$

\bigskip

The maximum correlation in (\ref{maxcorr}) remains valid for all cases 
$m \ge d \ge 0$.  Further, the maximum correlation converges to $1$ as 
$n \to \infty$, so we obtain $\mathrm{Var}(U-V)=(1-\mathrm{corr}(U,V))/6 \to 0$,
or $V-U\to 0$ in probability.  Therefore, as $n\to\infty$, the random 
variable $(U,V)=(U,U+(V-U))$ converges in law to $(U,U)$, a bivariate 
random variable whose joint distribution, remarkably, happens to provide
the Fr\'{e}chet--Hoeffding upper bound, $\min\{u,v\}$.
Thus, we have the following result.

\begin{theorem}
\label{thm:FHbound}
Let $C^*$ be the maximum correlation copula function (\ref{C*}) that 
is constructed by the B-spline 
\[
\{N^d_{k-d-1}\}_{k=1}^n = \{N^d_{i}: i=-d,-d+1,\ldots,m-2,m-1\}
\]
on $[0,1]$ of degree $d\ge 0$,
having equally-spaced knots (\ref{knots}) with $t_i=i/m$, $i=0,1,\ldots,m$, 
where $m \ge d$.  
As $m\to\infty$, $C^*(u,v) \to \min\{u,v\}$ for all $u, v$,
the Fr\'{e}chet--Hoeffding upper bound.
\end{theorem}

Table \ref{tab:rankcorr} shows the maximum correlations when the number
of basis functions is $n$.  In view of Table \ref{tab:rankcorr},
the range of correlation for the B-spline copulas of small order $d$
is wider than that of the Bernstein copula.  Indeed, 
$$
\mathrm{corr}(U,V)\approx 1-\frac{d+1}{n^2}.
$$
On the other hand, $d$ determines the smoothness of the copula density.
Consequently, some criterion is needed to evaluate data fitness so
as to balance the accuracy of the approximation with the smoothness
of the density; this problem will be studied in future work.

\medskip

\begin{table}[ht]
\caption{Maximum correlations}
\label{tab:rankcorr}
\begin{center}
\begin{tabular}{c|ccccc}
\hline
& & & & & \\[-5pt]
       & Bernstein$^*$ & $d=0$ & $d=1$ & $d=2$ & $d=3$ \\ [3pt]
\hline
& & & & & \\[-5pt]
$n=2$  & 0.333 & 0.75  & 0.333 & NA    & NA \\ [2pt]
$n=3$  & 0.5   & 0.889 & 0.667 & 0.5$^*$    & NA \\ [2pt]
$n=4$  & 0.6   & 0.938 & 0.827 & 0.688 & 0.6$^*$ \\ [2pt]
$n=5$  & 0.667 & 0.96  & 0.896 & 0.796 & ~0.72$^{**}$ \\ [2pt]
$n=6$  & 0.714 & 0.972 & 0.931 & 0.867 & 0.796 \\ [2pt]
$n=7$  & 0.75  & 0.980 & 0.951 & 0.908 & 0.851 \\ [2pt]
$n=8$  & 0.778 & 0.984 & 0.963 & 0.933 & 0.892 \\ [2pt]
$n=9$  & 0.8   & 0.988 & 0.971 & 0.949 & 0.919 \\ [2pt]
$n=10$ & 0.818 & 0.99  & 0.977 & 0.960 & 0.937 \\ [2pt]
\hline
& & & & & \\[-5pt]
$n$    &
$1-\dfrac{2}{n+1}$ & $1-\dfrac{1}{n^2}$ & $1-\dfrac{2(3n-5)}{3(n-1)^3}$ & $1-\dfrac{6n-19}{2(n-2)^3}$ & $1-\dfrac{2(50n-231)}{25(n-3)^3}$ \\ [2pt]
& & & & & \\[-5pt]
\hline
\end{tabular} \\[1mm]
\footnotesize{*: Bernstein case ($m=n-d=1$), \ \ \ **: Example \ref{ex:(n,d)=(5,3)}.}
\end{center}
\end{table}

We conjecture that Theorem \ref{thm:FHbound} holds in more general settings.

\medskip

\begin{conjecture}
\label{conjecture_maxcorr}
Let $(U,V)$ be distributed as the maximum correlation distribution
(\ref{C*}) constructed by the B-spline 
\[
\{N^d_{k-d-1}\}_{k=1}^n = \{N^d_{i}: i=-d,-d+1,\ldots,m-2,m-1\}
\]
on $[0,1]$, of degree $d \ge 0$, 
with the knots (\ref{knots}). 
As $m\to\infty$ with $\max_{1\le i\le m}|t_{i-1}-t_{i}|\to 0$,
$\mathrm{corr}(U,V)$ converges to $1$; hence, for all $u, v$, $C^*(u,v)$ 
converges to $\min\{u,v\}$, the Fr\'{e}chet--Hoeffding upper bound.
\end{conjecture}

\section{The maximum correlation copula: Total positivity}
\label{sec:total_positivity}
\setcounter{equation}{0}

The next two results improve significantly the previous ones about
the Bernstein copulas.

\begin{theorem}
\label{thm:TP-Hn}
The copula $C^*$ in (\ref{C*}) is TP${}_{\infty}$, i.e., 
for any $r\ge 1$,
\[
 \det\big(C^*(u_i,v_j)\big)_{r\times r} \ge 0
\]
for all $u_1>\cdots>u_r$ and $v_1>\cdots>v_r$.
\end{theorem}

\noindent
\textit{Proof.}  
All determinants arising in the proof are of order $r$, unless
otherwise specified.  Further, we consider two cases: (I) $r>n$, 
and (II) $r\le n$.

Case I: $r>n$.  In this case, the $r\times r$ matrix 
$(C^*(u_i,v_j)\bigr)_{1\le i,j\le r}$ satisfies 
\begin{align*}
\bigl(C^*(u_i,v_j)\bigr)_{1\le i,j\le r}
 &= \biggl(\sum_{k=1}^{n} q_k\Phi_{k}(u_i) \Phi_{k}(v_j)\biggr)_{1\le i,j\le r}\\
 &=\bigl(q_j\Phi_j(u_i)\bigr)_{1\le i\le r; 1\le j\le n}\bigl(\Phi_i(v_j)\bigr)_{1\le i\le n; 1\le j\le r}.
\end{align*}
Consequently, the rank of this matrix is at most $n$, and hence is degenerate.  Therefore, it follows obviously that $\det\big(C^*(u_i,v_j)\big) = 0$.  

Case II: $r\le n$.  We will show that $\det\big(C^*(u_i,v_j)\big) \ge 0$.  By
the Binet--Cauchy formula, 
\begin{align}
\label{binet-cauchy}
\det\big(C^*(u_i,v_j)\big)
&=\det\biggl(\sum_{k=1}^{n} q_k \Phi_{k}(u_i) \Phi_{k}(v_j)\biggr) \nonumber \\
&=\mathop{\sum\cdots\sum}_{n\ge k_1>\cdots>k_r\ge 1} \biggl( \prod_{i=1}^r q_{k_i} \biggr)
\det\bigl(\Phi_{k_i}(u_j)\bigr) \det\bigl(\Phi_{k_i}(v_j)\bigr).
\end{align}
Writing
\begin{equation}
\label{detPhi}
 \det(\Phi_{k_i}(u_j)) = \det\biggl(\int^{u_j}_{0} \phi_{k_i}(t)\, \dd t\biggr)
= \det\left(\int^{1}_{0} \phi_{k_i}(t) \1_{(0,u_j)}(t)\, \dd t\right),
\end{equation}
it follows by the continuous version of the Binet--Cauchy formula
\cite{Gross-Richards98,Karlin68} that 
\begin{multline}
\label{cvBCF}
\det\biggl(\int^{1}_{0} \phi_{k_i}(t) \1_{(0,u_j)}(t)\, \dd t\biggr) \\
= \mathop{\int\cdots\int}_{1>t_1>\cdots>t_r>0}
 \det\bigl(\phi_{k_j}(t_i) \bigr) \det\bigl(\1_{(0,u_i)}(t_j)\bigr) \prod_{i=1}^{r}\,\dd t_i.
\end{multline}

By Theorem \ref{thm:WT}, $\{\phi_{k_1},\ldots,\phi_{k_r}\}$ is an 
OCWT-system, hence 
\begin{equation}
\label{wT'}
 \det\bigl( \phi_{k_j}(t_i) \bigr) \ge 0
\end{equation}
for all $k_1>\cdots>k_r$ and $t_1>\cdots>t_r$.  
  Also, it is well-known from \cite{Karlin68,Karlin-Studden66} that
\[
 \det\bigl(\1_{(0,u_i)}(t_j)\bigr)\ge 0
\]
for all $u_1>\cdots>u_r$ and $t_1>\cdots>t_r$.  

Therefore, we deduce from (\ref{detPhi}) and (\ref{cvBCF}) that 
$\det\big(\Phi_{k_i}(u_j)\big) \ge 0$ for all 
$k_1>\cdots>k_r$ and $u_1>\cdots>u_r$. Similarly, we obtain
$\det\big(\Phi_{k_i}(v_j)\big) \ge 0$ for $k_1>\cdots>k_r$ and
$v_1>\cdots>v_r$. Hence, it follows from (\ref{binet-cauchy}) that
$\det\big(C^*(u_i, v_j)\big) \ge 0$ for $u_1>\cdots>u_r$ and
$v_1>\cdots>v_r$. The proof is complete.
\hfill $\qed$

\begin{remark}
For the case of the Bernstein copula, we note that (\ref{wT'}) is 
proved as follows.  Consider
\[
 \phi_k(t) = b_{k,n}(t) = k \binom{n}{k} t^{k-1} (1-t)^{n-k},
\]
$t \in [0,1]$.  Then,
\[
 \det(\phi_{k_j}(t_i)) = \biggl(\prod_{i=1}^{r} k_i \binom{n}{k_i}\biggr)
 \det\bigl(t^{k_j-1}_{i} (1-t_{i})^{n-k_j}\bigr)
\]
and
\begin{align*}
\det\Bigl(t^{k_j-1}_{i} (1-t_{i})^{n-k_j}\Bigr)
&= \det\Bigl(\Bigl(\frac{t_i}{1-t_i}\Bigr)^{k_j-1} (1-t_i)^{n-1} \Bigr) \\
&= \prod_{i=1}^{r} (1-t_i)^{n-1} \cdot \det\Bigl(\Bigl(\frac{t_i}{1-t_i}\Bigr)^{k_j-1}\Bigr).
\end{align*}

For $k_1>\cdots>k_r\ge 1$, set $k_j-1=\kappa_j+r-j$, $j=1,\ldots,r$,
and define the {\it partition} $\kappa=(\kappa_1,\ldots,\kappa_r)$, i.e., 
$\kappa_1,\ldots,\kappa_r$ are nonnegative integers and 
$\kappa_1 \ge \cdots \ge \kappa_r$.  
Also, let $z_i=t_i/(1-t_i)$, $i=1,\ldots,r$, and let
$z=(z_1,\ldots,z_r)$.  Recall from \cite{Macdonald15} that the {\it
Schur function} corresponding to the partition $\kappa$ is defined
as
$$
\chi_{\kappa}(z) = \frac{\det\big(z^{k_j-1}_{i}\big)}{\prod_{i<j} (z_i-z_j)}.
$$
Then we obtain
\begin{align*}
\det\big(z^{k_j-1}_{i}\big)
&= \chi_{\kappa}(z) \cdot \prod_{i<j} (z_i-z_j) \\
&= \chi_{\kappa}(z) \cdot \prod_{i<j} \Bigl(\frac{t_i}{1-t_i}-\frac{t_j}{1-t_j} \Bigr) \\
&= \chi_{\kappa}(z) \cdot \prod_{i<j} \frac{t_i-t_j}{(1-t_i)(1-t_j)}.
\end{align*}
It is well-known that $\chi_{\kappa}(z)\ge 0$ for 
and $z_1,\ldots,z_r \ge 0$ \cite{Gross-Richards98,Macdonald15}, and hence
\[
 \det\big(t^{k_j-1}_{i}(1-t_{i})^{n-k_j}\big) \ge 0
\]
for all $k_1>\cdots>k_r$ and $t_1>\cdots>t_r$.  This completes the proof
of (\ref{wT'}). 
\hfill $\qed$
\end{remark}

\medskip

As a consequence of Theorem \ref{thm:TP-Hn}, we obtain a new proof of 
the total positivity of the function $\min\{u,v\}$; see \cite[Chapter 2]{Karlin68}.

\begin{corollary}
\label{corFrechetHotelling}
The Fr\'{e}chet--Hoeffding upper bound, $\min\{u,v\}$, is TP${}_{\infty}$.
\end{corollary}

\noindent
\textit{Proof.}  
Recall that  for the Bernstein copula, $C^*(u,v)$ increases to
$\min\{u,v\}$ as the number $n$ of basis functions goes to infinity
(see Huang, et al. \cite{Huang-etal13}).  
Moreover, it follows from Theorem \ref{thm:FHbound} that for the
equally-spaced knot B-spline copula, $C^*(u,v)$ converges to
$\min\{u,v\}$.  In either case, by taking the limit, as $n\to\infty$, of
the $r \times r$ nonnegative determinant, $\det\big(C^*(u_i,v_j)\big)$, 
we obtain 
$$
\min\{u,v\} = \lim_{n \to \infty} C^*(u,v) \ge 0,
$$
which proves that the function $\min\{u,v\}$ is TP$_r$.  
Finally, since $r$ is arbitrary then it follows that the function 
$\min\{u,v\}$ is TP$_\infty$.
\hfill $\qed$

\medskip

By mimicking the proof of Theorem \ref{thm:TP-Hn}, we actually have the
following stronger result.  The proof is omitted.

\begin{theorem}
\label{thm:TP-c*}
The copula density $c^*$ in (\ref{c*}) is TP${}_{\infty}\,$.
\end{theorem}

Theorem \ref{thm:TP-Hn} is in fact a consequence of Theorem
\ref{thm:TP-c*} by using Lemma \ref{lem:ldk} below, but we provide a
direct proof there.

Let $(X,Y)\sim H$ with marginal distributions $F$ and $G$, and copula 
function $C$. Using the language of reliability theory, define the 
\textit{survival functions} 
$$
\overline{F}(x) = 1-F(x), \qquad 
\overline{G}(y) = 1-G(y),
$$
and 
$$
\overline{H}(x,y)=\Pr(X>x,\,Y>y)=1-F(x)-G(y)+H(x,y),
$$
$x, y \in \mathbb{R}$.  It follows from the definition of the copula
function that $H(x,y) = C(F(x),G(y))$ and
$\overline{H}(x,y)=\overline{C}(\overline{F}(x),\overline{G}(y))$
for all $x,y\in\mathbb{R}$. Recently, Lin et al. \cite{Lin-etal18}
proved the following result.

\begin{lemma}
\label{lem:ldk}
If the bivariate distribution $H$ has a TP$_r$ density with $r\ge 2$,
then both $H$ and $\overline{H}$ are TP$_r\,$.
Consequently, if $H$ has a TP$_{\infty}$ density,
both $H$ and $\overline{H}$ are TP$_{\infty}$.
\end{lemma}

An immediate consequence of the last two theorems is the following
result.  In part \textit{(ii)} of this result, we apply the fact that both $F$ 
and $G$ are non-decreasing, while both $\overline{F},\ \overline{G}$ 
are non-increasing (see Marshall, et al. 
\cite[p. 758]{Marshall-Olkin-Arnold11}).

\begin{corollary}
Let $C^*$ be the copula defined in (\ref{C*}). \\
(i)
The survival function $\overline{C^*}$ is
TP${}_{\infty}\,$. \\
(ii) 
If $(X,Y)\sim H$ with copula $C^*$, then both $H$ and
$\overline{H}$ are TP${}_{\infty}\,$.
\end{corollary}

We next discuss some implications of the total positivity.  By the
results of Gross and Richards \cite[Section 3, Example 3.7]{Gross-Richards98} we have
the following inequalities.

\begin{corollary} Let $(X,Y)\sim H$ with marginals $F, G$ and  copula $C^*$
in (\ref{C*}) and $r\ge 2$. Then the matrix
\begin{equation}
\label{EXiYjmatrix}
 \bigl( E[X^{i-1}Y^{j-1}] \bigr)_{1\le i,j\le r} =
 \begin{pmatrix}
 1          & E[Y]         & \cdots & E[Y^{r-1}] \\
 E[X]       & E[X Y]       & \cdots & E[X Y^{r-1}] \\
 \vdots     & \vdots       &        & \vdots \\
 E[X^{r-1}] & E[X^{r-1} Y] & \cdots & E[X^{r-1} Y^{r-1}]
\end{pmatrix}
\end{equation}
is TP$_r$, provided the expectations exist.

Let $x_1<\cdots<x_r$ and $y_1<\cdots<y_r$.  The matrix
\begin{equation}
\label{Hxiyjmatrix}
 \bigl( \overline{H}(x_i,y_j) \bigr)_{1\le i,j\le r} =
 \begin{pmatrix}
 \overline{H}(x_1,y_1) & \cdots & \overline{H}(x_1,y_r) \\
 \vdots      &        & \vdots \\
 \overline{H}(x_r,y_1) & \cdots & \overline{H}(x_r,y_r)
\end{pmatrix}
\end{equation}
is TP$_r$.
\end{corollary}

\medskip

In particular, when $r=3$, it follows from (\ref{EXiYjmatrix}) that 
\[
 \det\begin{pmatrix}
 1      & E[Y]     & E[Y^2] \\
 E[X]   & E[X Y]   & E[X Y^2] \\
 E[X^2] & E[X^2 Y] & E[X^2 Y^2]
\end{pmatrix}
\ge 0,
\]
an inequality that is equivalent to 
\begin{align*}
E[&X^2 Y^2] E[X Y] - E[X^2 Y] E[X Y^2] - E[X] E[X^2 Y^2] E[Y] \\
& + E[X^2] E[X Y^2] E[Y] + E[X] E[X^2 Y] E[Y^2] - E[X^2] E[XY] E[Y^2] \ge 0.
\end{align*}

Let $x_1=-\infty<x_2=x<x_3=x'$ and $y_1=-\infty<y_2=y<y_3=y'$.  Note
that $\overline{F}(x)=\overline{H}(x,-\infty)$ and
$\overline{G}(y)=\overline{H}(-\infty,y)$.  By (\ref{Hxiyjmatrix}), the 
matrix 
\[
 \begin{pmatrix}
 1      & \overline{G}(y)    & \overline{G}(y') \\
 \overline{F}(x)  & \overline{H}(x,y)  & \overline{H}(x,y') \\
 \overline{F}(x') & \overline{H}(x',y) & \overline{H}(x',y')
 \end{pmatrix}
\]
is totally positive of order $3$.  By calculating the $2 \times 2$ principal 
minor of this matrix, we find that 
$\overline{H}(x,y) \ge \overline{F}(x) \overline{G}(y)$; equivalently, 
$
 {H}(x,y) \ge {F}(x) {G}(y)$,
$x,y \in \mathbb{R}$, i.e., the distribution function ${H}$ is 
positively quadrant dependent.  Further, by calculating the determinant 
of this matrix, we obtain 
\begin{align*}
\overline{H}(&x',y')\overline{H}(x,y) - \overline{H}(x',y) \overline{H}(x,y') - \overline{F}(x) \overline{H}(x',y') \overline{G}(y) \\
& + \overline{F}(x') \overline{H}(x,y') \overline{G}(y) + \overline{F}(x) \overline{H}(x',y) \overline{G}(y') - \overline{F}(x') \overline{H}(x,y) \overline{G}(y') \ge 0
\end{align*}
for $x<x'$,\, $y<y'$.

We remark that more general inequalities can be deduced from \cite[Example 3.11]{Gross-Richards98}.

\section{Moments of the B-spline functions with initial boundary}
\label{sec:moment}
\setcounter{equation}{0}

 In this section, we provide the moment formula for the B-spline
functions with initial boundary at $t=0$ defined on
$\mathbb{R}_+=[0,\infty)$. The expressions for $q_k$ and $r_k$ in
Lemma \ref{lem:m} are obtained, in Corollary \ref{cor:h=0,1} below, 
as a consequence of the moment formula.  

Let $N^d_{i}$ be a B-spline function of degree $d \ge 0$ on 
$\mathbb{R}_+$ with knots:
\begin{equation}
\label{g-knots}
 \underbrace{t_{-d}=\cdots =t_{-1}}_d=t_0=0< t_1=1 < t_2=2 <\ \cdots
\end{equation}
(compare with the previously studied B-spline function 
defined in (\ref{knots})).  Here, we have $t_i = (i)_+ = \max\{i, 0\}$ 
and, as before, $N^d_{i}(t)$ is generated by the following recursion formula: 
\begin{equation}
\label{Ndi_recursion}
 N^d_i(t) = \frac{t-(i)_+}{(i+d)_+ -(i)_+}N^{d-1}_i(t) + \frac{(i+d+1)_+ -t}{(i+d+1)_+ -(i+1)_+}N^{d-1}_{i+1}(t),
\end{equation}
$d \ge 1$, with initial conditions
\[
 N^0_i(t) = \begin{cases} 1, & i\ge 0 \ \mbox{and}\ t\in [i,i+1), \\
 0, & \mbox{otherwise}. \end{cases}
\]
For each $i\ge -d$, $N^d_{i}$ is a non-zero function with support
$[\max\{i, 0\}, i+d+1]$.  The recurrence (\ref{Ndi_recursion}) can be 
written more concretely as 
\begin{align*}
 N^d_i(t)
& = \begin{cases}
\displaystyle
 \frac{t-i}{d}N^{d-1}_i(t) + \frac{i+d+1 -t}{d}N^{d-1}_{i+1}(t), & i\ge 0, \\[2mm]
\displaystyle
 \frac{t}{i+d}N^{d-1}_i(t) + \frac{i+d+1 -t}{i+d+1}N^{d-1}_{i+1}(t), & 
 -d<i\le -1, \\[2mm]
\displaystyle
 (1-t)N^{d-1}_{-d+1}(t), & i=-d, \\[1mm]
\displaystyle
 0, & i<-d.
\end{cases}
\end{align*}

For $h\ge 0$, denote the $h$-moment of $N^d_i$, 
\[
 \gamma^d_i(h) := \int_{-\infty}^\infty t^h N^d_i(t)\,\dd t = \int_{\max\{i,0\}}^{i+d+1} t^h N^d_i(t)\,\dd t;
\]
this quantity was used in the proof of Lemma \ref{lem:m} above.  
Then, we have the following recurrence relation for these moments.
\begin{equation}
\label{m-recur}
\gamma^d_i(h)
 =\begin{cases}
\displaystyle
 \frac{\gamma^{d-1}_i(h+1)-i \gamma^{d-1}_i(h)}{d} & \\[1mm]
\displaystyle\quad
 +\frac{(i+d+1)\gamma^{d-1}_{i+1}(h) -\gamma^{d-1}_{i+1}(h+1)}{d}, & i\ge 0, \\[2mm]
\displaystyle
 \frac{\gamma^{d-1}_i(h+1)}{i+d} & \\
\displaystyle\quad
+ \frac{(i+d+1)\gamma^{d-1}_{i+1}(h) -\gamma^{d-1}_{i+1}(h+1)}{i+d+1}, & 
-d<i<0, \\
\displaystyle
 \gamma^{d-1}_{-d+1}(h) -\gamma^{d-1}_{-d+1}(h+1), & i=-d, \\[1mm]
\displaystyle
 0, & i<-d,
\end{cases}
\end{equation}
with boundary condition
\begin{equation}
\label{boundary}
 \gamma^0_i(h) = \begin{cases} \displaystyle \frac{(i+1)^{h+1}-i^{h+1}}{h+1}, & i\ge 0, \\ 0, & i<0. \end{cases}
\end{equation}

The next result, which is interesting in its own right, presents
the solution of the recurrence system in terms of the Stirling 
numbers of the second kind:
\[
 S(n,k) = \frac{1}{k!}\sum_{j=0}^k (-1)^j \binom{k}{j} (k-j)^n.
\]
Here, $S(n,0)=\delta_{n0}$, $S(n,k)=0$ for $n<k$, and  $0^0\equiv 1$
whenever it arises.  Note also that $S(n,1) = S(n,n) = 1$ and 
$S(n,n-1) = n(n-1)/2$.  The Stirling numbers of the second kind 
satisfy the recurrence formula
\begin{equation}
\label{recurrence}
 S(n+1,k)=kS(n,k)+S(n,k-1),
\end{equation}
and the identity
\begin{equation}
\label{identity}
 S(n+1,k+1) = \sum_{j=k}^{n} \binom{n}{j}S(j,k) 
= \sum_{j=0}^{n-k} \binom{n}{n-j}S(n-j,k),
\end{equation}
which will be used later. For the identity (\ref{identity}), see
Wagner \cite{Wagner96} and the end of Remark 2 below.

\begin{theorem}
\label{thm:m-identity}
For $d \ge 0$, the $h$-moment of the B-spline
function $N^d_i$ in (\ref{g-knots}) is of the form
\begin{equation}
\label{m-identity}
\gamma^d_i(h) =
 \begin{cases}
\displaystyle
 \sum_{\ell=0}^h i^\ell \binom{h}{\ell}
 \dfrac{S(h+d+1-\ell,d+1)}{\displaystyle\binom{h+d+1-\ell}{d+1}}, & i\ge 0, \\[30pt]
 \displaystyle
 \dfrac{i+d+1}{d+1}\,\frac{S(h+i+d+1,i+d+1)}{\displaystyle\binom{h+d+1}{d+1}}, & -d\le i\le 0, \\[30pt]
 0, & i<-d.
\end{cases}
\end{equation}
\end{theorem}

\medskip

\begin{corollary}
\label{cor:h=0,1}
For $h=0,1$, we have
\begin{equation}
\label{h0}
\gamma^d_i(0) = \begin{cases}
 1, & i\ge 0, \\[10pt]
\displaystyle
 \dfrac{i+d+1}{d+1}, & -d\le i\le 0, \\[10pt]
 0, & i<-d
\end{cases}
\end{equation}
and
\begin{align}
\gamma^d_i(1)
& = \begin{cases}
\displaystyle
 \dfrac{d+2i+1}{2}, & i\ge 0, \\[10pt]
\displaystyle
 \dfrac{(i+d+1)^2(i+d+2)}{2(d+1)(d+2)}, & -d\le i\le 0, \\[10pt]
 0, & i<-d.
\end{cases}
\label{h1}
\end{align}
\end{corollary}

\medskip

The formulas (\ref{h0}) and (\ref{h1}) can be applied to obtain the formula
for the maximum correlation (\ref{maxcorr}).

\begin{corollary}
\label{cor:i=0,1}
For $i=0,1$,
\begin{equation}
\label{i=0,1}
 \gamma^d_i(h) = \dfrac{S(h+i+d+1,i+d+1)}{\displaystyle\binom{h+d+1}{d+1}}.
\end{equation}
\end{corollary}

\noindent
\textit{Proof.}  
For the case $i=1$, 
\begin{align*}
\gamma^d_1(h)
&= \sum_{\ell=0}^h \binom{h}{\ell} \dfrac{S(h+d+1-\ell,d+1)}{\displaystyle\binom{h+d+1-\ell}{d+1}} \\
&= \dfrac{1}{\displaystyle\binom{h+d+1}{d+1}} \sum_{\ell=0}^h \binom{h+d+1}{h+d+1-\ell} S(h+d+1-\ell,d+1) \\
&= \dfrac{S(h+d+2,d+2)}{\displaystyle\binom{h+d+1}{d+1}},
\end{align*}
by the identity (\ref{identity}).
\hfill $\qed$

\medskip

\noindent
\textit{Proof of Theorem \ref{thm:m-identity}}.  
We prove the statement by mathematical induction on $d$.
Note first that (\ref{m-identity}) with $d=0$ coincides with 
the boundary conditions (\ref{boundary}) for all $i$ and $h$.

Suppose that (\ref{m-identity}) is true for the case $d-1$ and
for all $i$ and $h$, 
then we wish to prove that it also holds true for
the case $d$ and for all $i$ and $h$.

(i) For $i \ge 0$, by the assumption of induction,
\[
 \gamma^{d-1}_{i}(h) = \sum_{\ell=0}^h i^\ell \binom{h}{\ell} \frac{S(h+d-\ell,d)}{\displaystyle\binom{h+d-\ell}{d}}.
\]
Then, we have
$$
\gamma^{d-1}_{i+1}(h) = \sum_{k=0}^h (i+1)^k \binom{h}{k} \frac{S(h+d-k,d)}{\displaystyle\binom{h+d-k}{d}},
$$
and by expanding $(i+1)^k$ using the binomial theorem, we obtain 
$$
\gamma^{d-1}_{i+1}(h) = \sum_{k=0}^h \sum_{\ell=0}^k i^\ell \binom{k}{\ell} \binom{h}{k} \frac{S(h+d-k,d)}{\displaystyle\binom{h+d-k}{d}}.
$$
Interchanging the order of summation and using the identity, 
$$
\binom{k}{\ell} \binom{h}{k} = \binom{h}{\ell} \binom{h-\ell}{k-\ell},
$$
we find that 
$$
 \gamma^{d-1}_{i+1}(h) = \sum_{\ell=0}^h i^\ell \binom{h}{\ell} \sum_{k=\ell}^h \binom{h-\ell}{k-\ell} \frac{S(h+d-k,d)}{\displaystyle\binom{h+d-k}{d}}.
$$
Replacing $k$ by $k-l$, we have 
$$
 \gamma^{d-1}_{i+1}(h) = \sum_{\ell=0}^h i^\ell \binom{h}{\ell} \sum_{k=0}^{h-\ell} \binom{h-\ell}{k} \frac{S(h+d-\ell-k,d)}{\displaystyle\binom{h+d-\ell-k}{d}},
$$
and using the identity, 
$$
\frac{\displaystyle\binom{h-\ell}{k}}{\displaystyle\binom{h+d-\ell-k}{d}} = \frac{\displaystyle\binom{h+d-\ell}{h+d-\ell-k}}{\displaystyle\binom{h+d-\ell}{d}},
$$
we deduce that 
\begin{align}
 \gamma^{d-1}_{i+1}(h) &= \sum_{\ell=0}^h i^\ell \binom{h}{\ell} \frac{1}{\displaystyle\binom{h+d-\ell}{d}}
\sum_{k=0}^{h-\ell} \binom{h+d-\ell}{h+d-\ell-k} S(h+d-\ell-k,d) \nonumber \\
&= \sum_{\ell=0}^h i^\ell \binom{h}{\ell} \frac{S(h+d-\ell+1,d+1)}{\displaystyle\binom{h+d-\ell}{d}},
\label{i+1}
\end{align}
where the last equality follows from the identity (\ref{identity}).  
Moreover, using the identity, 
$$
\binom{h+1}{\ell} = \binom{h}{\ell-1} + \binom{h}{\ell},
$$
we obtain 
\begin{align}
\label{gamma_identity}
\gamma^{d-1}_i&(h+1)-i \gamma^{d-1}_i(h) \nonumber \\
&= i^{h+1} + \sum_{\ell=0}^h \Biggl[ i^\ell \binom{h+1}{\ell} 
\frac{S(h+d-\ell+1,d)}{\displaystyle\binom{h+d-\ell+1}{d}} \nonumber \\
& \qquad\qquad\qquad\qquad\qquad\qquad\qquad\quad 
- i^{\ell+1} \binom{h}{\ell} \frac{S(h+d-\ell,d)}{\displaystyle\binom{h+d-\ell}{d}} \Biggr] \nonumber \\
&= i^{h+1} + \sum_{\ell=0}^h \Biggl[ i^\ell \biggl\{ \binom{h}{\ell-1} + \binom{h}{\ell} \biggr\} 
\frac{S(h+d-\ell+1,d)}{\displaystyle\binom{h+d-\ell+1}{d}} \nonumber \\
& \qquad\qquad\qquad\qquad\qquad\qquad\qquad\quad 
- i^{\ell+1} \binom{h}{\ell} \frac{S(h+d-\ell,d)}{\displaystyle\binom{h+d-\ell}{d}} \Biggr].
\end{align}
Since $S(d,d) = 1$ then 
\begin{align*}
i^{h+1} + \sum_{\ell=0}^h i^\ell \binom{h}{\ell-1}\frac{S(h+d-\ell+1,d)}
{\displaystyle\binom{h+d-\ell+1}{d}} 
&= \sum_{\ell=1}^{h+1} i^\ell \binom{h}{\ell-1}\frac{S(h+d-\ell+1,d)}
{\displaystyle\binom{h+d-\ell+1}{d}} \\
&= \sum_{\ell=0}^{h} i^\ell \binom{h}{\ell}\frac{S(h+d-\ell,d)}
{\displaystyle\binom{h+d-\ell}{d}},
\end{align*}
and substituting this result into (\ref{gamma_identity}), we obtain 
\begin{equation}
\label{gamma_identity2}
\gamma^{d-1}_i(h+1)-i \gamma^{d-1}_i(h) = \sum_{\ell=0}^h i^\ell \binom{h}{\ell} \frac{S(h+d-\ell+1,d)}{\displaystyle\binom{h+d-\ell+1}{d}}.
\end{equation}

Similarly, it follows from (\ref{i+1}) that
\begin{align}
\label{gamma_identity3}
\gamma^{d-1}_{i+1}&(h+1)-i\gamma^{d-1}_{i+1}(h) \nonumber \\
&= i^{h+1} + \sum_{\ell=0}^h \Biggl[ i^\ell \binom{h+1}{\ell} \frac{S(h+d-\ell+2,d+1)}{\displaystyle\binom{h+d-\ell+1}{d}} \nonumber \\
& \qquad\qquad\qquad\qquad\qquad\qquad
 - i^{\ell+1} \binom{h}{\ell} \frac{S(h+d-\ell+1,d+1)}{\displaystyle\binom{h+d-\ell}{d}} \Biggr] \nonumber\\
&= \sum_{\ell=0}^h i^\ell \binom{h}{\ell} \frac{S(h+d-\ell+2,d+1)}{\displaystyle\binom{h+d-\ell+1}{d}}.
\end{align}
Hence, by substituting (\ref{gamma_identity2}) and (\ref{gamma_identity3}) into (\ref{m-recur}), we find that 
\begin{align*}
\gamma^{d}_{i}(h)
&= \frac{\gamma^{d-1}_i(h+1)-i \gamma^{d-1}_i(h)}{d}
 -\frac{\gamma^{d-1}_{i+1}(h+1)-i\gamma^{d-1}_{i+1}(h)}{d}
 +\frac{d+1}{d}\gamma^{d-1}_{i+1}(h) \\
&= \sum_{\ell=0}^h i^\ell \binom{h}{\ell} \Biggl[ \frac{S(h+d-\ell+1,d)}{d \displaystyle\binom{h+d-\ell+1}{d}} - \frac{S(h+d-\ell+2,d+1)}{d \displaystyle\binom{h+d-\ell+1}{d}} \\
& \qquad\qquad\qquad\qquad\qquad\qquad
 + \frac{(d+1)S(h+d-\ell+1,d+1)}{d \displaystyle\binom{h+d-\ell}{d}} \Biggr] \\
&= \sum_{\ell=0}^h i^\ell \binom{h}{\ell} \frac{S(h+d-\ell+1,d+1)}{\displaystyle\binom{h+d-\ell+1}{d+1}},
\end{align*}
where the last equality follows from the recurrence formula
(\ref{recurrence}).

\smallskip

(ii) When $-d<i<0$, by the assumption of induction,
\[
 \gamma^{d-1}_{i}(h) = \frac{i+d}{d} \cdot \frac{S(h+i+d,i+d)}{\displaystyle\binom{h+d}{d}}.
\]
It then follows from (\ref{m-recur}) that
\begin{align*}
\gamma^{d}_{i}(h)
=& \frac{\gamma^{d-1}_i(h+1)}{i+d} +
\frac{(i+d+1)\gamma^{d-1}_{i+1}(h) -\gamma^{d-1}_{i+1}(h+1)}{i+d+1} \\
=& \frac{1}{i+d} \frac{i+d}{d}\,\frac{S(h+i+d+1,i+d)}{\displaystyle\binom{h+d+1}{d}} \\
&\qquad + \frac{i+d+1}{d}\,\frac{S(h+i+d+1,i+d+1)}{\displaystyle\binom{h+d}{d}} \\
&\qquad -\frac{1}{i+d+1} \frac{i+d+1}{d}\,\frac{S(h+i+d+2,i+d+1)}{\displaystyle\binom{h+d+1}{d}} \\
=& \frac{i+d+1}{d+1}\,\frac{S(h+i+d+1,i+d+1)}{\displaystyle\binom{h+d+1}{d+1}}.
\end{align*}
Here we used the recurrence formula (\ref{recurrence}) with $n=h+i+d+1$ and $k=i+d+1$, viz., 
\[
 S(h+i+d+2,i+d+1) = (i+d+1)S(h+i+d+1,i+d+1) + S(h+i+d+1,i+d).
\]

(iii) When $i=-d$, by the inductive hypothesis, 
\[
 \gamma^{d-1}_{-d+1}(h) = \frac{h! (d-1)!}{(h+d)!}.
\]
Then, by (\ref{m-recur}), we have 
\begin{align}
\gamma^{d}_{-d}(h)
&= \gamma^{d-1}_{-d+1}(h) -\gamma^{d-1}_{-d+1}(h+1) 
 = \frac{h! \, d!}{(h+d+1)!},
\label{m:-d}
\end{align}
which coincides with (\ref{m-identity}) with $i=-d$.

The proof is completed by induction on $d$.
\hfill $\qed$

\medskip

\begin{remark}
Recall the generalized (higher-order) Bernoulli
polynomial $B_{\ell}^{(x)}$ defined by the generating function
\begin{equation}
\label{gbp}
 \biggl(\frac{t}{e^t-1}\biggr)^x=\sum_{\ell=0}^{\infty}B_{\ell}^{(x)}\frac{t^{\ell}}{{\ell}!},
\end{equation}
$|t|<2\pi$, $x\in\mathbb{R}$, where $B_{\ell}^{(x)}$ 
is a polynomial of degree $\ell$ in $x$ with
rational coefficients. Neuman \cite[Proposition 3.5]{Neuman81} 
showed that $\gamma^d_0(h)=B_h^{(-(d+1))}$, from which (\ref{i=0,1})
for $i=0$ in Corollary \ref{cor:i=0,1} also follows by the
relationship between the Stirling number and the generalized
Bernoulli polynomial:
\begin{eqnarray}
\label{Sidentity0}
 S(n+k,k) = \binom{n+k}{k} \, B_{n}^{(-k)}.
\end{eqnarray}
The latter can be verified by (\ref{gbp}) and the exponential 
generating function, 
\begin{eqnarray}
\label{Sidentity}
 \frac{1}{k!}(e^t-1)^k=\sum_{n=k}^{\infty}S(n,k)\frac{t^n}{n!},
\end{eqnarray}
$t\in{\mathbb R}$, $k\ge 0$.  See \cite{Carlitz60} and \cite{Branson00} 
for (\ref{Sidentity0}) and
 (\ref{Sidentity}), respectively.  The generating function (\ref{Sidentity}) 
 can also be established by induction on $k$ and this is equivalent 
to verifying the above useful identity (\ref{identity}).
\end{remark}


\begin{remark}
By iteration, we have 
\[
 N_{-d}^d(t)=(1-t)N_{-d+1}^{d-1}(t)=\cdots=(1-t)^dN_0^0(t)=(1-t)^d\1_{[0,1)}(t),
\]
and hence its $h$-moment is equal to
\[
\gamma_{-d}^d(h)
 =\int_0^1 t^h(1-t)^d\,\dd t 
 =\frac{\Gamma(h+1)\,\Gamma(d+1)}{\Gamma(h+d+2)}
 =\frac{h!\,d!}{(h+d+1)!},
\]
as shown in (\ref{m:-d}). 
\end{remark}

\begin{remark}
It can be shown that for $i\ge 0$,
\[
 N_{i+1}^d(t)=N_{i}^d(t-1),\quad i+1\le t<i+d+2.
\]
Therefore, the $h$-moment of $N_{i+1}^d$ can be calculated as
\begin{align*}
\gamma_{i+1}^d(h)
&=\int_{i+1}^{i+d+2}t^hN_{i+1}^d(t)\,\dd t=\int_{i+1}^{i+d+2}t^hN_{i}^d(t-1)\,\dd t \nonumber \\
&=\int_i^{i+d+1}(t+1)^hN_i^d(t)\,\dd t=\sum_{j=0}^h \binom{h}{j} \gamma_i^d(j),\ \ i\ge 0.
\end{align*}
This is equivalent to the first formula of (\ref{m-identity})
(with $i\ge 0$).  Indeed, for $i\ge 1$,
it follows from (\ref{m-identity}) that
\begin{align*}
 \sum_{j=0}^h \binom{h}{j} \gamma_{i-1}^d(j)
&= \sum_{j=0}^h \binom{h}{j} 
 \sum_{k=0}^j (i-1)^k  \binom{j}{k} 
 \frac{S(j+d+1-k,d+1)}{\displaystyle\binom{j+d+1-k}{d+1}},
 \end{align*}
and if we now change variables from $0\le k\le j\le h$ to 
$0\le k\le \ell\le h$, where $\ell=h-j+k$, then we obtain 
\begin{align*}
 \sum_{j=0}^h \binom{h}{j} \gamma_{i-1}^d(j)
&= \sum_{\ell=0}^h \sum_{k=0}^\ell
 (i-1)^k \binom{\ell}{k} \binom{h}{\ell}
 \frac{S(h+d+1-\ell,d+1)}{\displaystyle\binom{h+d+1-\ell}{d+1}} \\
&= \sum_{\ell=0}^h
 i^{\ell} \binom{h}{\ell}
 \frac{S(h+d+1-\ell,d+1)}{\displaystyle\binom{h+d+1-\ell}{d+1}} \\
 &= \gamma_{i}^d(h).
\end{align*}
\end{remark}

\subsection*{Acknowledgments}
This work was supported by JSPS KAKENHI Grants, Numbers 16K00060 and 16H02792.

\end{document}